\documentclass[12pt]{amsart}
\usepackage{amscd,amsmath,amsthm,amssymb}
\usepackage{pstcol,pst-plot,pst-3d}
\usepackage{lmodern,pst-node}
\def\NZQ{\Bbb}               
\def\NN{{\NZQ N}}

\def\ZZ{{\NZQ Z}}

\def\PP{{\NZQ P}}

\def\GG{{\NZQ G}}

\def\frk{\frak}               

\def\mm{{\frk m}}

\def\Phi{{\frk n}}
\def\Phi{{\frk N}}
\def\MP{{\mathcal P}}
\def\MQ{{\mathcal Q}}
\def\MR{{\mathcal R}}
\def\MV{{\mathcal V}}
\def\xb{{\bold x}}

\def\opn#1#2{\def#1{\operatorname{#2}}} 
\opn\chara{char} \opn\length{\ell} \opn\pd{pd} \opn\rk{rk}
\opn\projdim{proj\,dim} \opn\injdim{inj\,dim} \opn\rank{rank}
\opn\depth{depth} \opn\grade{grade} \opn\height{height}
\opn\embdim{emb\,dim} \opn\codim{codim}

\opn\Tr{Tr} \opn\bigrank{big\,rank}
\opn\superheight{superheight}\opn\lcm{lcm}
\opn\trdeg{tr\,deg}
\opn\reg{reg} \opn\lreg{lreg} \opn\ini{in} \opn\lpd{lpd}
\opn\size{size}\opn\bigsize{bigsize}
\opn\cosize{cosize}\opn\bigcosize{bigcosize}
\opn\sdepth{sdepth}\opn\sreg{sreg}
\opn\link{link}\opn\fdepth{fdepth}
\opn\div{div} \opn\Div{Div} \opn\cl{cl} \opn\Cl{Cl}
\opn\Spec{Spec} \opn\Supp{Supp} \opn\supp{supp} \opn\Sing{Sing}
\opn\Ass{Ass} \opn\Min{Min}\opn\Mon{Mon}
\opn\Ann{Ann} \opn\Rad{Rad} \opn\Soc{Soc}
\opn\Im{Im} \opn\Ker{Ker} \opn\Coker{Coker} \opn\Am{Am}
\opn\Hom{Hom} \opn\Tor{Tor} \opn\Ext{Ext} \opn\End{End}
\opn\Aut{Aut} \opn\id{id}

\opn\nat{nat}
\opn\pff{pf}
\opn\Pf{Pf} \opn\GL{GL} \opn\SL{SL} \opn\mod{mod} \opn\ord{ord}
\opn\Gin{Gin} \opn\Hilb{Hilb}\opn\sort{sort}
\opn\aff{aff} \opn\con{conv} \opn\relint{relint} \opn\st{st}
\opn\lk{lk} \opn\cn{cn} \opn\core{core} \opn\vol{vol}
\opn\link{link} \opn\star{star}\opn\lex{lex}
\opn\gr{gr}

\def\pot#1#2{#1[\kern-0.28ex[#2]\kern-0.28ex]}

\opn\dirlim{\underrightarrow{\lim}}
\opn\inivlim{\underleftarrow{\lim}}
\let\union=\cup
\let\sect=\cap
\let\dirsum=\oplus
\let\tensor=\otimes

\let\Sect=\bigcap
\let\Dirsum=\bigoplus

\let\to=\rightarrow
\let\To=\longrightarrow
\def\Implies{\ifmmode\Longrightarrow \else
        \unskip${}\Longrightarrow{}$\ignorespaces\fi}
\def\implies{\ifmmode\Rightarrow \else
        \unskip${}\Rightarrow{}$\ignorespaces\fi}
\def\iff{\ifmmode\Longleftrightarrow \else
        \unskip${}\Longleftrightarrow{}$\ignorespaces\fi}

\let\:=\colon
\newtheorem{Theorem}{Theorem}[section]
\newtheorem{Lemma}[Theorem]{Lemma}
\newtheorem{Corollary}[Theorem]{Corollary}
\newtheorem{Proposition}[Theorem]{Proposition}

\newtheorem{Example}[Theorem]{Example}
\newtheorem{Examples}[Theorem]{Examples}

\let\epsilon\varepsilon
\let\kappa=\varkappa
\def\qed{\ifhmode\textqed\fi
      \ifmmode\ifinner\quad\qedsymbol\else\dispqed\fi\fi}
\def\textqed{\unskip\nobreak\penalty50
       \hskip2em\hbox{}\nobreak\hfil\qedsymbol
       \parfillskip=0pt \finalhyphendemerits=0}
\def\dispqed{\rlap{\qquad\qedsymbol}}

\opn\dis{dis}
\def\pnt{{\raise0.5mm\hbox{\large\bf.}}}

\opn\Lex{Lex}

\begin{document}

\title {Stanley depth and size  of a monomial ideal}
\author {J\"urgen Herzog, Dorin  Popescu and Marius Vladoiu}
\address{J\"urgen Herzog, Fachbereich Mathematik, Universit\"at Duisburg-Essen, Campus Essen, 45117
Essen, Germany} \email{juergen.herzog@uni-essen.de}

\address{Dorin Popescu, Institute of Mathematics "Simion Stoilow" of the Romanian Academy, University of Bucharest,
P.O.Box 1-764, Bucharest 014700, Romania}\email{dorin.popescu@imar.ro}

\address{Marius Vladoiu, Faculty of Mathematics and Computer Science, University of Bucharest, Str.
 Academiei 14, Bucharest, RO-010014, Romania}\email{vladoiu@gta.math.unibuc.ro}

\thanks{This paper was partially written during the visit of the first author at the Institute of Mathematics ``Simion Stoilow" of the Romanian Academy supported by a BitDefender Invited Professor Scholarship, 2010. The second and the third author were partially supported by the CNCSIS grant PN II-542/2009, respectively CNCSIS grant TE$\_46$ nr. 83/2010 of Romanian Ministry of Education, Research and Innovation, and also want to express their gratitude to ASSMS of GC University Lahore for creating a very appropriate atmosphere for research work.}

\begin{abstract}
Lyubeznik introduced the concept of size of a monomial ideal and showed that  the size of a monomial ideal increased by $1$ is a lower bound for its depth. We show that the size  is also a lower bound for its Stanley depth. Applying Alexander duality we obtain upper bounds for the regularity and Stanley regularity of squarefree monomial ideals.
\end{abstract}
\maketitle

\section*{Introduction}
Let $I\subset S$ be a monomial ideal where $S=K[x_1,\ldots,x_n]$ is  the polynomial ring  in the indeterminates $x_1,\ldots,x_n$ over the field $K$.  In \cite{L} Lyubeznik showed that $\depth I\geq  1+\size I$. In the case that $I$ is a squarefree monomial ideal with  minimal prime ideals  $P_1,\ldots,P_s$, the size of $I$ is the number $v+(n-h)-1$, where $h$ is the height of $\sum_{j=1}^sP_j$ and $v$ is the minimal number $t$ for which there exist integers $i_1<i_2<\cdots <i_t$ such that $\sum_{k=1}^tP_{i_k}=\sum_{j=1}^sP_j$. Replacing in the previous definition ``there exist integers $i_1<i_2<\cdots <i_t$'' by  ``for all  integers $i_1<i_2<\cdots <i_t$'', one obtains the definition of bigsize, first considered by Popescu in \cite{Dorin}.  See Section 1 where the general definition of size and bigsize is given.

It is conjectured by Stanley that $\sdepth I\geq \depth I$. Assuming Stanley's conjecture and combining it with Lyubeznik's inequality one should expect that  $\sdepth I\geq 1+ \size I$, where $\sdepth I$ denotes the Stanley depth of $I$. In Section 3 we show that this is indeed the case, see Theorem~\ref{sizeineq}.  We also expect that
$\sdepth S/I\geq \size I$. Assuming the conjectured inequality $\sdepth I\geq 1+\sdepth S/I$, our inequality $\sdepth I\geq 1+ \size I$ would also follow from the  inequality $\sdepth S/I\geq \size I$.

In Section 1 we introduce a cocomplex $\GG$ which is a attached to a set ${\mathcal P}=\{J_1,\ldots,J_s\}$ of monomial ideals in $S/I$ where $I\subset S$ is a monomial ideal. This cocomplex is acyclic (Theorem~\ref{acyclic}), and in the case that $I=0$ and the ideals $J_k$ are irreducible monomial ideals, $\GG$ may be viewed as the Alexander dual of the Taylor complex. The complex allows in some cases to compute or to estimate the depth of $\Sect_{k=1}^sJ_k$. In terms of this complex we give in Theorem~\ref{sizedepth} a criterion for a monomial ideal $I$ to have minimal depth, that is, to satisfy the equation $\depth I=1+\size I$, which is in particular the case when $\bigsize I=\size I$. We conclude Section~1 by giving upper bounds  for the regularity of a squarefree monomial ideal $I$ in terms of the cosize of $I$. This upper bound is obtained by applying Alexander duality, see Corollary~\ref{dual}.

In Section 2 we describe the method  of splitting the variables in order to obtain lower bounds for the Stanley depth of monomial ideals.  This technique was first  introduced  by A.~Popescu \cite{Adi} in a special case, and then generalized by D.~Popescu \cite{Dorin} for all squarefree monomial ideals. Here we further extend it to arbitrary monomial ideals, and use it in Section 3 to prove the inequality $\sdepth I\geq 1+ \size I$. We conclude the paper by giving upper bounds for the Stanley regularity in terms of the cosize of $I$.

\section{Size, big size and depth of a squarefree monomial ideal}

Let $K$ be a field and $S=K[x_1,\ldots,x_n]$ be the  polynomial ring over $K$ in the indeterminates $x_1,\ldots,x_n$. Let $I\subset S$ be a squarefree monomial ideal, and $I=\Sect_{i=1}^sP_i$ its presentation as an irredundant intersection of irreducible  monomial ideals. This unique presentation establishes  a bijection between  monomial  ideals and finite sets of irreducible monomial ideals, that is, monomial ideals of the form $(x_{i_1}^{a_1},x_{i_2}^{a_2},\ldots,x_{i_k}^{a_k})$,   with no inclusions among them. By means of this intersection we would like to give some depth estimates for $I$.

More generally, let $\MP=\{I_1,\ldots, I_s\}$ be an arbitrary set of  monomial ideals, and $I\subset S$ a monomial ideal. We set $T=S/I$.  Given these data, we are going to introduce a (co)complex $\GG$ of $T$-modules which allows in some cases to compute the depth of $J=(\Sect_{j=1}^s I_j+I)/I$. For $k=1,\ldots,s$, let $J_k=(I_k+I)/I\subset T$. Then the  $i$th component $G_i$ of $\GG$ is given by
\[
G_i=\Dirsum_{U\subset [s], |U|=i}T_U,
\]
where $T_U=T/(\sum_{k\in U}J_k)$. The chain map from $G_i\to G_{i+1}$  is defined on the component $T_U\to T_V$ to be the canonical residue class map multiplied with $(-1)^\ell$, were $\ell$ is the number of elements $i\in U$ with $i<k$ if $V=U\union \{k\}$,  and to be the zero map if $U\not\subset V$.

For $s=2$ and $I=0$ the complex  $\GG$ is the standard complex
\[
\begin{CD}
0@>>>  S @>>>  S/I_1\dirsum S/I_2 @>>> S/(I_1+I_2)  @>>>  0\\
&&a+I_1\sect I_2& \mapsto & (a+I_1,a+I_2)\\
&&&&(a+I_1,b+I_2) & \mapsto & (a-b)+I_1+I_2.
\end{CD}
\]
In general, $\GG$ is of the form
\[
0\to T\to \Dirsum_{i=1}^s T_i\to \Dirsum_{i<j}T_{ij}\to \cdots \to T_{12\cdots s}\to 0.
\]
It is easy to see  that $\GG$ is indeed a complex. We call $\GG$ the complex attached to the set $\MP=\{J_1,\ldots, J_s\}$ of monomial ideals in $T=S/I$. Now we have
\begin{Theorem}
\label{acyclic}
Let $\GG$ be the complex attached to $\MP=\{J_1,\ldots,J_s\}$. Then
\[
H_i(\GG)= \left\{ \begin{array}{ll}
       J, & \;\text{if  $i=0$}, \\ 0, & \;\text{otherwise,}
        \end{array} \right.
        \]
        where $J=\Sect_{k=1}^sJ_k$.
\end{Theorem}

\begin{proof}
We prove the theorem by induction on $s$. The assertion is trivial for $s=1$. Now let $s>1$, and let $\GG'$  be the complex attached to the set $\MQ=\{J_1,\ldots,J_{s-1}\}$.  Furthermore, let  $\GG''$ be  the complex attached to the set of monomial ideals $$\MR=\{(J_1+J_s)/J_s,\ldots,(J_{s-1}+J_s)/J_s\}\quad  \text{in}\quad  T/J_s.$$ Then $\GG''$ may be viewed as a subcomplex of $\GG$, and we get a short exact sequence of complexes
\[
0\To  \GG''[-1]\To \GG\to\GG'\To 0,
\]
where $\GG''[-1]$ is the complex $\GG''$ homologically shifted by $-1$.

This short exact sequence gives rise to the  following long exact sequence
\[
\cdots \To H_{i-1}(\GG'')\To H_{i}(\GG)\To H_{i}(\GG')\To \cdots.
\]
By our induction hypothesis we have $H_{i}(\GG')=H_{i}(\GG'')=0$ for $i>0$. This implies that $H_i(\GG)=0$ for $i>1$, and in addition we get the exact sequence
\[
0\To H_{0}(\GG)\To H_{0}(\GG')\To H_{0}(\GG'')\To H_{1}(\GG)\To 0.
\]
Our induction hypothesis implies that  $H_{0}(\GG')=\Sect_{k=1}^{s-1}J_k$ and that  $$H_{0}(\GG'')=\Sect_{k=1}^{s-1}(J_k+J_s)/J_s.$$ The map $H_{0}(\GG') \to H_{0}(\GG'')$ is  just the canonical  map $$\Sect_{k=1}^{s-1}J_k\to \Sect_{k=1}^{s-1}(J_k+J_s)/J_s.$$ Since the ideals $J_k$ are all monomial ideals in $T=S/I$ it follows that $\Sect_{k=1}^{s-1}(J_k+J_s)=(\Sect_{k=1}^{s-1}J_k)+J_s$. This implies that $H_{0}(\GG'') \to H_{0}(\GG')$ is surjective, so that $H_1(\GG)=0$, and that $\Ker(H_{0}(\GG'') \to H_{0}(\GG'))=J$, so that $H_0(\GG)=J$, as desired.
\end{proof}

Let $I\subset S$  be a monomial ideal and $I=\Sect_{i=1}^sQ_i$ a primary decomposition
of $I$, where the $Q_i$ are monomial ideals. Let $Q_i$ be $P_i$-primary. Then each $P_i$ is a monomial prime ideal and $\Ass(S/I)=\{P_1,\ldots,P_s\}$.

According to  Lyubeznik \cite[Proposition 2]{L} the {\em size} of $I$, denoted $\size I$,
 is the number $v + (n - h) - 1$,  where $v$ is  the minimum number $t$ such that there exist $j_1<\cdots< j_t$  with
\[
\sqrt{\sum_{k=1}^tQ_{j_k}}=	\sqrt{\sum_{j=1}^sQ_j},
\]
and where $h=\height \sum_{j=1}^sQ_j$.

Notice that $\sqrt{\sum_{k=1}^tQ_{j_k}} =\sum_{k=1}^tP_{j_k}$ and $\sqrt{\sum_{j=1}^sQ_j}=\sum_{j=1}^sP_j$, so that the size of $I$ depends only on the set of associated prime ideals of $S/I$.

If in the above definition of $\size I$ we replace ``there exist $j_1<\cdots< j_t$" by ``for all $j_1<\cdots< j_t$", we obtain the definition of $\bigsize I$, introduced by Popescu \cite{Dorin}. Of course $\bigsize I\geq \size I$, and in fact the big size of $I$ is in general much bigger than the size of $I$.

To illuminate these concepts consider the so-called squarefree Veronese ideal $I_{n,d}\subset S$ consisting of all squarefree monomials of degree $d$. It is known that the minimal prime ideals of $I_{n,d}$ are exactly all the monomial prime ideals of height $n-d+1$ and that $\depth I_{n,d}=d$. From this information one easily deduces that $\size I_{n,d}=\lceil n/(n-d+1)\rceil$ and $\bigsize I_{n,d}=d-1$. Thus we see that
\[
\bigsize I_{n,d}+1=\depth I_{n,d}\geq \size I_{n,d}+1.
\]
In general $\depth I$ can be bigger or smaller than $\bigsize I+1$. The inequality  $\depth I_{n,d}\geq \size I_{n,d}+1$ is just a special case of Lyubeznik's inequality \cite[Proposition 2]{L}.
We say that $I$ has {\em minimal depth}, if equality holds, i.e. $\depth I=\size I+1$. In the last  section  it  will be shown that for a  monomial ideal one has  $\sdepth I\geq 1+\size I$. Thus if $I$  has  minimal depth, then  $\sdepth I\geq \depth I$, in which case Stanley's conjecture holds.

In the  next result  a sufficient condition is given for a monomial ideal to achieve the lower bound for the depth, as  given by Lyubeznik. Recall that a monomial ideal  is irreducible if and only if it is generated by powers of subsets of the variables. Each monomial $I$ has a unique presentation $I=\Sect_{j=1}^sQ_j$ as an intersection of minimal irreducible monomial ideals. Moreover, $\Ass(I)=\{P_1,\ldots,P_s\}$, where $P_i=\sqrt{Q_i}$. For example,  for $I=(x^2,xy,x^2z^2,xyz^2,y^2z^2)$ we have $I= (x^2,y)\sect (x,y^2)\sect (x,z^2)$ and $\Ass(I)=\{(x,y),(x,z)\}$.

For the proof of the next result it is important to notice that any sum of irreducible monomial ideals is again irreducible, and hence a complete intersection. In particular, for any subset $U\subset [s]$ we have $\dim S_U=\depth S_U$ where as before, $S_U=S/(\sum_{j\in U}Q_j)$.

\begin{Theorem}
\label{sizedepth}
Let $I\subset S$ be a monomial ideal of size $t$ and $\MP=\{Q_1,\ldots,Q_s\}$ be the  set of minimal irreducible prime ideals of $I$, and assume without loss of generality that $\sum_{j=1}^sQ_j$ is $\mm$-primary, where $\mm$ is the graded maximal ideal of $S$.

Let $\GG$ be the complex attached to $\MP$ in $S$.  Let $\MV=\{V\: \; |V|=t+1,\; \dim S_V=0\}$,  and let $W=\Dirsum_{V\subset \MV}S_V$. Assume that
\begin{enumerate}
\item[(a)] $\Im(G_t\to G_{t+1})\sect W\neq (0)$;
\item[(b)] for all $U\subset [s]$ with $|U|=t$ there exists $V\subset \MV$ with $U\subset V$.
\end{enumerate}
 Then $I$ has minimal depth. In particular, the conditions of the theorem are satisfied if $\bigsize I=\size I$.
\end{Theorem}

\begin{proof}
We may assume that $\sqrt{\sum_{j=1}^s Q_j}=\mm$, where $\mm$ is the graded maximal ideal of $S$. Let $U_{i+1}=\Im(G_i\to G_{i+1})$. We are going to show that $\depth U_i= t+1-i$. Then, since $U_1=S/I$, the desired conclusion follows. The module   $U_{t+1}$ has a non-trivial intersection with $W$. Since each element of $W$ is annihilated by a power of  $\mm$, it follows  from (a) that $U_{t+1}$ contains a nonzero element which is annihilated by $\mm$. This shows that $\depth U_{t+1}=0$.

Next we will show by induction on $i$ that $\depth G_{t+1-i}\geq i$. (The desired formula for the depth of the $U_i$ is then obviously a consequence of this fact.) For $i=1$ the assertion follows from the assumption that $\size I=t$. Let us now assume  that $\depth G_{t+1-i}\geq i$, but $\depth G_{t+1-(i+1)}<i+1$. Then there exists $U\subset [s]$ with $|U|=t-i$ and $\dim S_U\leq i$. It follows from condition (b)  that there exists $V\subset [s]$ with $|V|=t+1$ such that $U\subset V$ and $\dim S_V=0$. Choose such a set $V$ and an element $k\in V\setminus U$. Since  $\depth G_{t+1-i}\geq i$ it follows that $\dim S_{U\union\{k\}}\geq i$. On the other hand, it is clear that $\dim S_U\geq \dim S_{U\union\{k\}}$. Hence we conclude that $\dim S_U=\dim S_{U\union\{k\}}$ which is only possible if $Q_k\subset \sqrt{\sum_{j\in U}Q_j}$. It follows that  $\sqrt{\sum_{j\in V}Q_j}=\sqrt{\sum_{j\in V\setminus \{k\}}Q_j}$, contradicting the fact that $\size I=t$.
\end{proof}

\begin{Examples}
{\em (i) Let $I=(x_1,x_2,x_3)\sect (x_1,x_4)\sect (x_2,x_4)\sect (x_3,x_4)$. Then $\depth I=1+\size I=2$, and hence $I$ has minimal depth. In fact, $I$ satisfies the hypotheses of Theorem~\ref{sizedepth}. On the other hand, $\bigsize I=2>\size I$.

(ii) Let $\Delta$ be the simplicial complex on the vertex set $\{1,\ldots,6\}$, associated to the canonical  triangulation of the real projective plane $\PP^2$, whose facets are
\[
\mathcal F(\Delta)=\{125,126,134,136,145,234,235,246,356,456\}.
\]
Then the Stanley-Reisner ideal of $\Delta$ (see definition below) is
\[
I_{\Delta}=(x_1x_2x_3,x_1x_2x_4,x_1x_3x_5,x_1x_4x_6,x_1x_5x_6,x_2x_3x_6,x_2x_4x_5,x_2x_5x_6,x_3x_4x_5,x_3x_4x_6).
\]
It is known that $\depth I_{\Delta}=4$ if $\chara K\neq 2$ and $\depth I_{\Delta}=3$ if $\chara K=2$. On the other hand, it is easy to see that $\size I_\Delta=2$ and $\bigsize I_\Delta=3$. This shows that the property of a monomial ideal to have minimal depth may depend on the characteristic of the base field. It can also be easily checked that condition (b) of Theorem~\ref{sizedepth} is satisfied for the ideal $I_\Delta$. On the other hand, since $I$ does not have minimal depth  if $\chara K\neq 2$, it follows that in this case condition (a) is not satisfied for $I_\Delta$. In particular, (b) does not imply (a).

(iii) Let $I=(x_1,x_2,x_3)\sect(x_2,x_3,x_4)\sect(x_1,x_4,x_5,x_6)\sect(x_1,x_3,x_5,x_7)\sect(x_2,x_4,x_6,x_7)$. Then $\depth I=3$ and $\size I=1$. The ideal satisfies condition (a) but not condition (b). Thus the examples (ii) and (iii) show that the condition (a) and (b) in Theorem~\ref{sizedepth} are independent.
}
\end{Examples}

By using Alexander duality one easily obtains a statement which is dual to that of Lyubeznik and also dual to that of Theorem~\ref{sizedepth}.

In order to describe it, let, as before, $S = K[x_1, \ldots, x_n]$ be the polynomial ring
in $n$ variables over a field $K$ and $\Delta$ a simplicial
complex on $[n]$.  For each subset $F \subset [n]$ we set
\[
\xb_F = \prod_{i \in F} x_i.
\]

Recall that the Stanley--Reisner ideal of $\Delta$ is the ideal
$I_\Delta$ of $S$ which is generated by those squarefree monomials
$\xb_F$ with $F \not\in \Delta$. One sets $K[\Delta]=S/I_\Delta$.

Then the Alexander dual of $\Delta$ is defined to be the simplicial complex
\[
\Delta^\vee = \{ [n] \setminus F : F \not\in \Delta \}.
\]
Obviously one has $(\Delta^\vee)^\vee=\Delta$.

We quote the following fact which for example can be found in \cite{HH}.
\begin{itemize}
\item
Let $I_{\Delta} = P_{F_1} \cap \cdots
\cap P_{F_m}$ be the standard primary decomposition of
$I_\Delta$. (Here $P_G=(\{x_i\}_i\in G)$ for $G\subset [n])$. Then  $\{ x_{F_1}, \ldots, x_{F_m} \}$ is the minimal monomial set of generators of $I_{\Delta^\vee}$.

\item (Terai) $\projdim I_\Delta=\reg K[\Delta^\vee]$, where as usual $\reg M$ denotes the regularity of a finitely generated graded $S$-module $M$.
\end{itemize}

Let $I$ be a squarefree monomial ideal minimally generated by the monomials $u_1,\ldots,u_m$. Let  $w$ be the smallest number $t$ with the property that  there exist integers $1\leq i_1<i_2<\cdots<i_{t}\leq m$ such that
\[
\lcm(u_{i_1},u_{i_2},\ldots,u_{i_{t}})=\lcm(u_{1},u_{2},\ldots,u_{m}).
\]
Then we call the number $\deg \lcm(u_{1},u_{2},\ldots,u_{m})- w$ the {\em cosize} of $I$,  denoted $\cosize I$. If in the above definition we replace the  words  `there exist'  by `for all', then we obtain the definition of the {\em big cosize} of $I$, denoted $\bigcosize I$.

Now we have
\begin{Corollary}
\label{dual}
Let $I\subset S$ be a squarefree monomial ideal. Then
\begin{enumerate}
\item[(a)] $\reg S/I\leq \cosize I$.
\item[(b)] $\reg S/I=\cosize I$, if $\bigcosize I=\cosize I$.
\end{enumerate}
\end{Corollary}

\begin{proof}
Let $\Delta$ be the simplicial complex with the property that $I=I_\Delta$. By using the result of Lyubeznik as well as the above facts, we obtain
\begin{eqnarray*}
n-\reg K[\Delta]= n-\projdim I_{\Delta^\vee}=\depth I_{\Delta^\vee}\geq \size I_{\Delta^\vee}+1,
\end{eqnarray*}
so that $n-\reg K[\Delta]\geq v+(n-h)$. This implies that $\reg K[\Delta]\leq h-v$.

Since $(\Delta^\vee)^\vee=\Delta$, we see that the number $v$ for $I_{\Delta^\vee}$ is equal to the number $w$ for $I_\Delta$, and that the number $h$ for $I_{\Delta^\vee}$ is equal to the number $\deg\lcm(u_{1},u_{2},\ldots,u_{m})$ for $I_\Delta$. Thus statement (a) follows.

The assertion (b) is a simple consequence of (a) and Theorem~\ref{sizedepth}.
\end{proof}

\section{Splitting the variables to get lower bounds for the Stanley depth}

In this section we describe and  extend a method, introduced in the papers \cite{Adi} and \cite{Dorin}, to decompose a  monomial ideal $I\subset S$ into $\ZZ^n$-graded  subspaces which allows us  to bound from below the Stanley depth of a monomial ideal. The decomposition depends on the choice of a subset $Y$ of the set of variables $X=\{x_1,\ldots,x_n\}$, and is also determined by the unique irredundant presentation of $I$ as an intersection $I=\Sect_{j=1}^sQ_j$  of its minimal irreducible monomial ideals. As before each $Q_j$ is a $P_j$-primary ideal.

Without loss of generality we may assume that $Y=\{x_1,\ldots,x_r\}$ for some number  $r$ such that $0\leq r\leq n$. Then the set of variables splits into the two sets $\{x_1,\ldots,x_r\}$ and $\{x_{r+1},\ldots,x_n\}$.

Given a subset $\tau\subset [s]$, we let $I_\tau$ be the $\ZZ^n$-graded $K$-vector space spanned by the  set of monomials of the form $w=uv$ where $u$ and $v$ are monomials with
\[
u\in K[x_1,\ldots,x_r]\quad \text{and}\quad  u\in \Sect_{j\not\in\tau}Q_j\setminus \sum_{j\in \tau}Q_j,
\]
\[
 v\in K[x_{r+1},\ldots,x_n]\quad \text{and}\quad v\in  \Sect_{j\in\tau}Q_j.
\]

The following result extends the corresponding statement shown by Popescu \cite{Dorin} for squarefree monomial ideals.

\begin{Proposition}
\label{newdorin}
With the notation introduced,   the ideal $I$ has a decomposition $\mathcal{D}_Y\: I =\Dirsum_{\tau\subset [s]}I_\tau$ as a direct sum  of $\ZZ^n$-graded $K$-subspaces of $I$.
\end{Proposition}

\begin{proof}
It is clear from the definition of $I_\tau$ that $I_\tau\subset I$, so that $\sum_{\tau\subset [s]}I_\tau\subset I$. Conversely, let $w=x_1^{a_1}x_2^{a_2}\cdots x_n^{a_n}$ be a monomial in $I$. Then $w=x_1^{a_1}x_2^{a_2}\cdots x_n^{a_n}$ can be written in a unique way as a product $w=uv$ of monomials   with $u\in K[x_1,\ldots,x_r]$ and $v\in K[x_{r+1},\ldots,x_n]$. Let $\tau=\{j\in[s]\: \; u\not \in Q_j\}$. Then $u\in \Sect_{j\not\in\tau}Q_j$ and $u\not \in \sum_{j\in \tau}Q_j$.

Let $j\in \tau$. Since $uv\in I$, it follows that $uv\in Q_j$. Thus, if  $Q_j=(x_{i_1}^{b_{i_1}},\ldots,x_{i_k}^{b_{i_k}})$, then there exists an integer $\ell$ with $a_{i_\ell}\geq b_{i_\ell}$. On the other hand, since $u\not\in Q_j$, it follows that $a_{i_t}<b_{i_t}$ for all $i_t\leq r$. This implies that $i_\ell\geq r+1$, and consequently $v\in Q_j$. Hence we see that $v\in \Sect_{j\in\tau}Q_j$, and conclude that $w\in I_\tau$.

In order to see that the sum is direct assume that $w=uv\in I_\tau\sect I_\sigma$. Then $u\in K[x_1,\ldots,x_r]$ and $u\in \Sect_{j\not \in \tau}Q_j\sect\Sect_{j\not\in\sigma}  Q_j=\Sect_{j\not\in \tau\union\sigma}Q_j$. Suppose that $\tau\neq \sigma$. Then we may assume that  $\sigma\setminus\tau\neq \emptyset$. Let $j\in \sigma\setminus\tau$. Then $u\in Q_j$ be the definition of $I_\tau$, and $u\not\in Q_j$, by the definition of $I_\sigma$, a contradiction.
\end{proof}

The $\ZZ^n$-graded $K$-subspaces  $I_\tau$ of $I$ have the structure of a $\ZZ^n$-graded module over $S$, and can be interpreted as follows: let
\[
S' =K[x_1,\ldots,x_r] \quad \text{and}\quad S''=K[x_{r+1},\ldots,x_n].
\]
Let $S\to S/(x_{r+1},\ldots,x_n)=S'$ be the canonical epimorphism and let $Q_i'$, $P_i'$ be the images of $Q_i$, respectively $P_i$ for
$i=1,\ldots,s$. Then we set
\[
J_\tau=\Sect_{j\not\in\tau} Q_j',
\]
and let $H_\tau$ be the $K$-vector subspace of $S'$ generated by all the monomials of $J_\tau\setminus(\sum_{j\in\tau}Q_j')$. Notice that $\sum_{j\in\tau}Q_j'$ is of the form $(x_{i_1}^{a_1},\cdots, x_{i_t}^{a_t})$ with $1\leq i_1<i_2<\cdots <i_t\leq r$ and suitable exponents $a_j>0$. Thus
\[
\label{idealfactors}
H_\tau=(J_\tau, x_{i_1}^{a_1},\cdots, x_{i_t}^{a_t})/(x_{i_1}^{a_1},\cdots, x_{i_t}^{a_t}).
\]

Next let $S\to S/(x_1,\ldots,x_r)=S''$ be the canonical epimorphism and let $Q_i''$ be the image of $Q_i$ for
$i=1,\ldots,s$. Then we define  the monomial ideal $L_\tau\subset S''$ as
\[
L_\tau=\Sect_{\j\in\tau}Q_j''.
\]
Now $I_\tau$ can be written as follows
\[
\label{tensor}
I_\tau=H_\tau\tensor_K L_\tau.
\]

The following example describes the decomposition of $I$ given in Proposition~\ref{newdorin}.

\begin{Example}
\label{tau-decomposition}
{\em Let $\mm$ be the maximal graded ideal of $K[x_1,x_2,x_3]$ and $\mm^2=Q_1\sect Q_2\sect Q_3$ be its unique irredundant irreducible decomposition. Then $Q_1=(x_1^2,x_2,x_3)$, $Q_2=(x_1,x_2^2,x_3)$ and $Q_3=(x_1,x_2,x_3^2)$. We choose the set $Y$ to be $\{x_1\}$, which splits the set of variables into the sets $\{x_1\}$ and $\{x_2,x_3\}$. Then $I=\Dirsum_{\tau\subset [3]} I_\tau$, where $I_\tau$ is the $K$-vector space generated by all monomials $w=uv$, where $u,v$ are monomials with
\[
u\in K[x_1] \quad \text{ and } \quad u\in\Sect_{j\not\in\tau}Q_j\setminus\sum_{j\in\tau} Q_j,
\]
\[
v\in K[x_2,x_3]\quad \text{ and } \quad v\in\Sect_{j\in\tau} Q_j.
\]
It follows immediately that the only nonzero summands $I_\tau$ of $I$ correspond to the following subsets of $\{1,2,3\}$:
\[
\tau\in\{\emptyset,\{1\},\{1,2,3\}\}.
\]
Indeed, since $Q_i+Q_j=\mm$ for all $i\neq j$ we obtain that $Q_k\setminus (Q_i+Q_j)=\emptyset$ for every permutation set $\{i,j,k\}$ of $\{1,2,3\}$ and therefore $I_\tau=0$ for all $\tau\in\{\{1,2\},\{1,3\},\{2,3\}\}$. For $\tau=\{3\}$ we obtain that $((Q_1\sect Q_2)\setminus Q_3)\sect K[x_1]=\emptyset$ and consequently $I_{\{3\}}=0$. Similarly, we obtain for $\tau=\{2\}$ that $I_{\{2\}}=0$. Now, for the remaining subsets we compute $I_\tau$. If $\tau=\emptyset$ then $I_{\emptyset}$ is generated as a $K$-vector space by all monomials $u\cdot v$, where
\[
u\in I\sect K[x_1]=(x_1^2)K[x_1] \quad \text{ and }\quad v\in K[x_2,x_3].
\]
Therefore $I_{\emptyset}=(x_1^2)K[x_1,x_2,x_3]$. For $\tau=\{1\}$, the $K$-basis of $I_{\{1\}}$ is given by the monomials $u\cdot v$, where
\[
u\in ((Q_2\sect Q_3)\setminus Q_1)\sect K[x_1]=x_1K \text{ and } v\in Q_1\sect K[x_2,x_3]=(x_2,x_3)K[x_2,x_3].
\]
Consequently we have that $I_{\{1\}}=(x_1x_2,x_1x_3)K[x_2,x_3]$. Finally, if $\tau=[3]$ we obtain that  $I_{[3]}$ is generated as  a $K$-vector space by all monomials $u\cdot v$, where
\[
u\in (S\setminus\mm)\sect K[x_1]=K \quad \text{ and }\quad v\in I\sect K[x_2,x_3]=(x_2^2,x_2x_3,x_3^2)K[x_2,x_3].
\]
Therefore $I_{[3]}=(x_2^2,x_2x_3,x_3^2)K[x_2,x_3]$. Hence we obtain the following decomposition of $I$ into $\ZZ^n$-graded $K$-subspaces of $I$:
\[
I=(x_1^2)K[x_1,x_2,x_3]\dirsum (x_1x_2,x_1x_3)K[x_2,x_3]\dirsum (x_2^2,x_2x_3,x_3^2)K[x_2,x_3].
\]
}
\end{Example}

\section{A comparison of Stanley depth and size}

As an application of the technique of splitting variables, as introduced in the previous section, we show

\begin{Theorem}
\label{sizeineq}
Let $I$ be a monomial ideal of $S$. Then  $$\sdepth I\geq 1+\size I.$$
\end{Theorem}
\begin{proof}
Let $I=\Sect_{j=1}^s Q_j$ be the unique irredundant presentation of $I$ as an intersection of its minimal irreducible monomial ideals. Each of the  $Q_j$ is a primary ideal whose associated monomial prime ideal we denote, as before, by $P_j$.

We may assume that $\sum_{j=1}^s P_j=\mm$. Indeed, let $Z=\{x_i\not\in  \sum_{j=1}^s P_j\}$,  $T=K[X\setminus Z]$ and $J=I\sect T$. Then the sum of the associated prime ideals of $J$ is the graded maximal ideal of $T$, and
\[
\sdepth I=\sdepth J +|Z|, \quad \text{and}\quad \size I=\size J +|Z|.
\]
The first equation follows from \cite[Lemma 3.6]{HVZ}, while the second equation follows from the definition of size.

 We choose the  splitting set $Y$ to be the set $\{x_i:x_i\in P_1\}$, and we may assume that $Y=\{x_1,\ldots,x_r\}$ for some number $r$ such that $1\leq r\leq n$. If $r=n$,  then the desired inequality follows at once since in this case $\size I=0$, and since for every monomial ideal $I$ we have that $\sdepth I\geq 1$. Therefore from now on we assume that $r<n$. We will prove the assertion of the theorem  by induction on $s$. The case $s=1$ follows immediately from \cite[Theorem 2.4]{Shen} and \cite[Lemma 3.6]{HVZ} since $\sdepth Q_1=\lceil |Y|/2\rceil + n-|Y|$ and $\size Q_1=n-|Y|$.

Assume now that the assertion is proved for all monomial ideals which are intersections of at most  $s-1$ irreducible monomial ideals.  Since $Y=\{x_1,\ldots,x_r\}$,  it follows from the method described before  Proposition~\ref{newdorin} that $I =\Dirsum_{\tau\subset [s]}I_\tau$ with $I_{[s]}=0$. We obtain from the decomposition of $I$ that
\[
\sdepth_S I\geq \min\{\sdepth_{S} I_\tau: \; \tau\subset [s]\; \text{and}\; I_\tau\neq 0 \}.
\]
Hence it remains to prove that for any subset $\tau$ of $[s]$ such that $I_\tau\neq 0$ we have that $\sdepth_{S} I_\tau\geq 1+\size I$. We will distinguish two cases:  $\tau=\emptyset$,  or $\tau$ is a proper non-empty subset of $[s]$. In both cases we may assume that $I_\tau\neq 0$.

In the first case we have that $I_\emptyset=(I\sect K[x_1,\ldots,x_r])S$. Applying now \cite{HVZ} and the fact that the sdepth of any ideal is greater than or equal to $1$ we obtain
\begin{eqnarray*}
\sdepth_{S} I_\emptyset&=&\sdepth_{K[x_1,\ldots,x_r]} (I\sect K[x_1,\ldots,x_r]) + n-r\geq 1 +\dim S/P_1\\
&\geq& 1+\depth S/I\geq 1+\size I.
\end{eqnarray*}

In the second case we first  observe that
\begin{eqnarray}
\label{sum}
\sdepth_{S} I\geq\min_\tau\{\sdepth_{S'} H_{\tau} + \sdepth_{S''} L_\tau\},
\end{eqnarray}
where we set $\sdepth M=0$ if $M=0$.

Indeed, we noticed already that $\sdepth_S I\geq \min\{\sdepth_S I_\tau\:\; \tau\subset [s]\}$. Since $I_\tau =H_\tau\tensor_K L_\tau$, it follows that $\sdepth I_\tau\geq \sdepth_{S'} H_\tau+\sdepth_{S''} L_\tau$, because if $\Dirsum_i u_iK[Z_i]$ is a Stanley decomposition of $H_\tau$, and $\Dirsum_j v_jK[W_j]$ is a Stanley decomposition of $L_\tau$, then $\Dirsum_{i,j}u_iK[Z_i]\tensor_K v_jK[W_j]$ is a Stanley decomposition of $I_\tau =H_\tau\tensor_K L_\tau$, see \cite[Lemma 1.2]{Adi}, where this assertion is shown in the case that $H_\tau$ and $L_\tau$ are both monomial ideals and \cite[Theorem 3.1]{Rauf} in the case  that $H_\tau$ and $L_\tau$ are both quotients of polynomial rings in disjoint sets of variables by monomial ideals. The argument in this slightly more general case is verbatim the same.

In our further discussions we distinguish whether  $P_1\not\subset \sum_{j\in\tau} P_j$ or $P_1\subset \sum_{j\in\tau} P_j$.

In the case that $P_1\not\subset \sum_{j\in\tau}P_j$, one may assume that
$H_\tau=J_\tau\setminus (x_k^{a_k},\cdots,x_r^{a_r})$ with $k>1$. In other words,
\begin{eqnarray}
\label{greater}
H_\tau =(J_\tau, x_k^{a_k},\cdots,x_r^{a_r})/(x_k^{a_k},\cdots,x_r^{a_r})
\end{eqnarray}
Thus $H_\tau$ is a submodule of $S'/(x_k^{a_k},\cdots,x_r^{a_r})$ from which it follows that $\depth H_\tau>0$. This in turn implies that $\sdepth H_\tau>0$, see \cite[Theorem 1.4]{C}. Since $|\tau|\leq s-1$,  $L_\tau$ is the intersection of at most $s-1$ irreducible monomial ideals. Thus, applying the induction hypothesis, the inequalities (\ref{sum}) and the subsequent  Lemma~\ref{smallerring}  we obtain
\[
\sdepth_{S} I_\tau\geq 1+\sdepth_{S''}L_\tau\geq 2+\size_{S''} L_\tau\geq 1+\size_S I,
\]
as desired.

On the other hand, if  $P_1\subset \sum_{j\in\tau} P_j$, then  $H_\tau$ has a presentation as in (\ref{greater}) but with $k=1$. Thus in this case $\depth H_\tau=0$. Again applying \cite[Theorem 1.4]{C} it follows that $\sdepth H_\tau=0$. Then as  before we get
\[
\sdepth_S I_\tau\geq \sdepth_{S''}L_{\tau}\geq 1+\size_{S''} L_\tau\geq 1+\size_S I,
\]
and we are done.
\end{proof}

\begin{Lemma}
\label{smallerring}
Let $I=\Sect_{i=1}^sQ_i$ be the unique irredundant presentation of $I$ as the intersection of irreducible monomial ideals, where the $Q_i$ are $P_i$-primary ideals. Assume that $P_1=(x_1,\ldots, x_r)$ is one of the minimal monomial prime ideals of $I$.  Let  $\tau\subset [s]$ such  that $L_{\tau}\neq 0$. Then
\[
\size_{S''}(L_\tau)+1\geq \size_S I,
\]
Moreover, if $P_1\subset\sum_{j\in\tau} P_j$ we even have
\[
\size_{S''}(L_\tau)\geq \size_S I.
\]
\end{Lemma}

\begin{proof}
We may assume as in Theorem~\ref{sizeineq} that $\sum_{j=1}^sP_j=\mm$. Let $c$ be the minimum number $t$ such that there exist $j_1<\cdots< j_t$ in $\tau$ with
\[
\sum_{k=1}^t (P_{j_k}\sect S'')=	\sum_{j\in\tau} (P_j\sect S'').
\]
We have to analyze two cases: $\sum_{j\in\tau}(P_{j}\sect S'')=(x_{r+1},\ldots,x_n)$ or $\sum_{j\in\tau}(P_{j}\sect S'')$ is properly contained in $(x_{r+1},\ldots,x_n)$. In the first case we have that $\size_{S''}(L_\tau)=c-1$ and $P_1+\sum_{k=1}^{c}P_{j_k}= P_1 + \sum_{j\in\tau}P_{j}=\mm$. This yields the first inequality. In particular, if $P_1\subset\sum_{j\in\tau} P_j$ then $\sum_{k=1}^{c}P_{j_k}= P_1 + \sum_{k=1}^{c}P_{j_k}=\mm$ and therefore we have the second inequality.

In the second case let $\{x_{i_1},\ldots,x_{i_d}\}$ be the variables from $S''$ that do not belong to $\sum_{j\in\tau} (P_j\sect S'')$. Then we have $\size_{S''}(L_\tau)=c-1+d$. Since, by our assumptions we have $P_1=(x_1,\ldots,x_r)$ and $\sum_{j=1}^s P_j=\mm$, it follows that for each $k$ with $1\leq k\leq d$ there exists an integer $l_k\in\{2,\ldots,s\}\setminus\tau$ such that $x_{i_k}\in P_{l_k}$. Then we have
\[
\sum_{k=1}^{c}(P_{j_k}\sect S'') + \sum_{k=1}^{d} (P_{l_k}\sect S'') = (x_{r+1},\ldots,x_n),
\]
and consequently
\[
P_1 + \sum_{k=1}^{c} P_{j_k} + \sum_{k=1}^{d} P_{l_k} = \mm.
\]
Hence we obtain that $\size I\leq c+d$, which is the desired first inequality. In particular, if $P_1\subset\sum_{j\in\tau} P_j$ then
\[
\sum_{k=1}^{c} P_{j_k} + \sum_{k=1}^{d} P_{l_k} = P_1 + \sum_{k=1}^{c} P_{j_k} + \sum_{k=1}^{d} P_{l_k} = \mm,
\]
which yields the second inequality.
\end{proof}

The reader may wonder why in the proof of Theorem~\ref{sizeineq} we have chosen the set $Y$ as the set of generators of one of the minimal prime ideals of $I$. This was chosen so to make sure that $I_{[s]}=0$. Indeed, if $I_{[s]}\neq 0$, then this would be a summand in the decomposition of $I$ which may have sdepth  less than or equal to the size of $I$, as the following  example  shows.

\begin{Example}
{\em Let $I=(x_1,x_2,x_3,x_6)\sect (x_2,x_3,x_4,x_6)\sect (x_2,x_3,x_5,x_6)$ be a monomial ideal of $K[x_1,\ldots,x_6]$. One can easily see that $\size  I=2$. If we choose now $Y$ to be the set $\{x_1,\ldots,x_5\}$ then the set of variables splits into the sets $\{x_1,\ldots,x_5\}$ and $\{x_6\}$. Then, for $\tau=[3]$ we have that $I_{[3]}$ is the $K$-vector space whose basis consists of the monomials $w=uv$, where $u,v$ are monomials with
\[
u\in (S \setminus\mm)\sect K[x_1,\ldots,x_5]=K \quad \text{ and } \quad  v\in I\sect K[x_6]=(x_6)K[x_6].
\]
Therefore $I_{[3]}=(x_6)K[x_6]$ and consequently $\sdepth I_{[3]}=1<\size I$.
}
\end{Example}

\medskip
\noindent
Dual to the  Lyubeznik inequality $\depth I\geq \size I +1$, we have $\reg I\leq \cosize I+1$, as  we have seen in Section 1. Similarly there is  an  inequality dual to $\sdepth I\geq \size I +1$, as we shall see now. For its proof we have to recall a few results.

\begin{enumerate}
\item[($\gamma$)]  Alexander duality  can be extended to finitely generated $\ZZ^n$-graded modules $M$, see \cite{Tim} and \cite{Y}. Then one  obtains a functor $M\mapsto M^\vee$ from the category of $\ZZ^n$-graded modules into itself with the property that $(M^\vee)^\vee=M$. Moreover one has  $(I_\Delta)^\vee=K[\Delta^\vee]$.
\item[($\delta$)] Let $\mathcal{D}\: M=\Dirsum_{i=1}^mu_iK[Z_i]$ be a Stanley decomposition of $M$. Then  \[\sreg \mathcal{D}=\max\{\deg u_i\:\; i=1,\ldots,m\}\] is called the Stanley regularity of $\mathcal{D}$, and
   \[
    \sreg M=\min\{\sreg \mathcal{D}\:\; \text{$\mathcal{D}$ is a Stanley decomposition of $M$}\}
    \]
is called the {\em Stanley regularity} of $M$. The crucial fact that we need has been shown by Soleyman Jahan \cite[Theorem 3.9]{Ali}, namely: $\sreg M=n-\sdepth M^\vee$.
\end{enumerate}

\begin{Corollary}
\label{sreg}
Let $I\subset S$ be a squarefree monomial ideal. Then  $\sreg S/I\leq \cosize I$, and equality holds if $\bigcosize I= \cosize I$.
\end{Corollary}

\begin{proof}
Let $\Delta$ be the simplicial complex with $I=I_\Delta$. Then, by using ($\gamma$) and ($\delta$) as well as Theorem~\ref{sizeineq}, we obtain
\[
\sreg S/I_\Delta=n-\sdepth I_{\Delta^\vee}\leq n- (\size I_{\Delta^\vee}+1)= \cosize  I_\Delta.
\]
\end{proof}


\begin{thebibliography}{10}
\bibitem{BH} W.\ Bruns and J. Herzog, {\em Cohen-Macaulay rings} Revised edition. Cambridge University Press (1998).


\bibitem{C} M.\ Cimpoea\c{s}, Some remarks on the Stanley depth for multigraded modules. Le Matematiche {\bf LXIII}, 165-175 (2008).



\bibitem{HH} J.\ Herzog and T.\ Hibi, {\em Monomial Ideals}. GTM 260. Springer 2010.

\bibitem{HVZ} J.\  Herzog,  M.\ Vladoiu, X.\  Zheng, How to compute the Stanley depth of a monomial ideal. J.~Algebra {\bf 322}, 3151--3169 (2009).

\bibitem{L} G.\ Lyubeznik, On the Arithmetical Rank of Monomial ideals. J.~Algebra {\bf 112}, 86--89  (1988).

\bibitem{Adi}  A.\ Popescu, Special Stanley Decompositions. Bull. Math. Soc. Sci. Math. Roumanie {\bf 53}(101), 361--372 (2010).

\bibitem{Dorin} D.\ Popescu, Stanley conjecture on intersections of four monomial prime ideals. arXiv.AC/1009.5646.

\bibitem{Rauf} A.\ Rauf, Depth and Stanley Depth of Multigraded Modules. Comm. in Algebra. {\bf 38}, 773--784 (2010).

\bibitem{Tim} T.\ R\"omer, Generalized Alexander duality and applications. Osaka J. Math.\  {\bf 38}, 469--485 (2001).

\bibitem{Shen} Y.\ Shen,  Stanley depth of complete intersection monomial ideals and upper-discrete partitions. J.~Algebra\ {\bf 321}, 1285--1292
(2009).

\bibitem{Ali} A.\ Soleyman--Jahan, Prime filtration and Stanley decompositions of squarefree modules and Alexander duality.
Manuscripta Math.\ {\bf 130}, 533--550 (2009).

\bibitem{Y} K.\ Yanagawa, Alexander duality for Stanley-Reisner rings and squarefree $\NN$-graded modules. J. Algebra,
{\bf 225},  630--645 (2000).

\end{thebibliography}
\end{document}